# An Extended AW-Rascle Model with Source Terms and Its Numerical Solution

Nandan Maiti1[0000-0002-2825-514X] and Bhargava Rama Chilukuri2[0000-0003-3772-796X]

1 LICIT-ECO7, University Gustave Eiffel, ENTPE, Lyon, France,69500
2 Department of Civil Engineering, IIT Madras, India, Chennai-600036
Email: nandan.besu043@gmail.com

**Abstract.** Nonlinear hyperbolic partial differential equations govern continuum traffic flow models. Higher-order traffic flow models consisting of continuum equations and velocity dynamics were introduced to address the limitations of the Lighthill, Whitham, and Richards (LWR) model. However, these models are ineffective in incorporating road heterogeneity. This paper integrates an extended AW-Rascle higher-order model with the source terms in the continuum equation to predict the traffic states in heterogeneous road conditions. The system of the equations was solved numerically with the central dispersion (CD) method incorporated into the standard McCormack scheme. Smoothing is applied to take care of the numerical oscillation of the higher-order model. Different combinations of initial conditions with source terms showed that the proposed model with the numerical methods could produce a stable solution and eliminate oscillation of the McCormack scheme.

## 1 Introduction

Macroscopic traffic flow models describe the dynamics and kinematics of aggregate traffic using mathematical relationships among the traffic flow variables in time-space. They commonly describe traffic flow as a continuum flow and are often compared to continuum models for compressible fluids. Continuum models assume traffic is similar to a compressible fluid flow. This family of models started with a well-known and oldest first-order time-space derivative LWR model (Lighthill & Whitham (1955); Richards (1956)) [14] and was later modified with higher-order models. The main standard equations for all these models are the conservation of vehicles and the fundamental equation of traffic flow. This model gives a dynamic equation for density in the time-space domain through (1). Equation 1 is a flow-conservative form of traffic, which is a nonlinear hyperbolic partial differential equation and can be rewritten as (2).

$$\frac{\partial k(x,t)}{\partial t} + \frac{\partial q(x,t)}{\partial x} = 0, \quad q = f(k) \tag{1}$$

$$\frac{\partial U}{\partial t} + A\frac{\partial U}{\partial x} = 0 \tag{2}$$

where density: $U = [k]$, flow: $F = [q]$, and wave-speed: $A = \frac{\partial F}{\partial U}$, where [.] represents the system of equations in matrix form.

The LWR model assumes single-lane homogeneous traffic. Several papers have attempted to describe heterogeneous traffic using two different approaches. The first approach was to make a distinction between lanes [22,21,9,11,12]. A second approach was to analyze the vehicle population based on their class, multi-class modeling [16,10,27]. Though the LWR model can explain the shockwave phenomena explicitly, it cannot explain other phenomena, such as stop-go waves, traffic hysteresis, phantom jams, etc.

A significant limitation of the LWR model is that it cannot capture the hysteresis effect. One of the main reasons for the hysteresis phenomena observed in traffic is the driver's response to the frontal stimuli and inertial effect. The higher-order models incorporated the inertial effects and driver anticipation on vehicle speed. Thus, these models are different from the LWR model, which has an additional equation of velocity dynamics. Payne [24] addressed one of the major drawbacks of the LWR model is that velocity adapts instantaneously to the traffic density by proposing a partial differential equation describing the velocity dynamics with the LWR system Eq. 1. Daganzo [4] reviewed Payne's model and stated that the kinematic theory of fluid should not be applied to model vehicular traffic. The main reason is traffic flow model should be anisotropic, where the fluid flow is isotropic. AW and Rascle (AWR) [2] proposed a new anisotropic continuum traffic flow model, addressing all the drawbacks of the Payne model [24] pointed out by Daganzo [4]. AWR introduced a convective derivative of pressure instead of a spatial derivative. This model considers not only the spatial changes of density but also the spatiotemporal changes in density.

All the higher-ordered models postulate the conservation of the number of vehicles and vehicle kinematics. According to fluid dynamics, a second conservation principle based on linear momentum ensures that the sum of all external forces acting on a system equals its time rate of linear momentum change. The flow of traffic seems not to be governed by any such principle. Also, in the case of traffic, link flow is adjusted by the entries and exits traffic from the main link. Therefore, the continuum equation must include the source terms in order to account for entries and exits and to correct the conservation of the velocity dynamics if unjustified.

The complex mathematical form of the higher-ordered traffic flow models calls for different numerical recipes to simplify it to a first-order model. Second-order continuum models in the presence of source terms are classified as inhomogeneous partial differential equations. McCormack method is a popular numerical framework for solving inhomogeneous continuum models due to its simplicity, limited computational costs, and less numerical difficulty [8]. Inspired by the hyperbolic second-order computational fluid dynamics and heat transfer application [13,1], a group of complex methods has also been proposed for continuum traffic flow applications [5,23]. The McCormack scheme is famous for having less numerical diffusion than a first-order model, e.g., Forward difference or backward difference [25]. However, this scheme is associated with numerical oscillation. Due to the lack of source terms, this scheme fails to reproduce feasible results for homogeneous higher-order models.

This research contributes to developing higher-order inhomogeneous traffic flow models by including source and sink terms in the continuum equations. The McCormack scheme has also been adjusted to reproduce rational results for different initial conditions when applied to an inhomogeneous second-order traffic flow model. Finally, the numerical test is validated with empirical data to check both the model's and numerical scheme's ability to predict traffic state for mixed traffic conditions. The main contributions of this work can be summarized as follows. We formulate a continuum equation of traffic flow with source terms and a velocity-dynamic equation in Section II. Then, we expressed the system of equations in a generalized form of PDE and proposed the numerical solutions for the PDE in section III. Section IV discussed the numerical test with different initial and source conditions. The proposed method is validated with the empirical observation in section V. Finally, Section VI presents discussions.

## 1.1    Motivation

Incorporating source terms in the continuum equation of traffic flow is motivated by the need to model the effects of external factors that can influence traffic dynamics. These external factors, which are not accounted for in the traditional traffic flow models, can include on-ramps, off-ramps, traffic signals, and other sources of congestion or delay. Incorporating source terms into the traffic flow equation allows for a more realistic and accurate representation of the traffic dynamics on the road network. For example, an on-ramp can introduce new vehicles into the flow, which can affect the density and velocity of the traffic stream. Similarly, a traffic signal can create stop-and-go waves that can propagate through the network, leading to congestion and delays.

Including source terms in the traffic flow equation can also improve the model's ability to predict the traffic system's behavior under different scenarios. For example, by incorporating the effects of a planned road closure or construction project, transportation planners can assess the potential impacts on traffic congestion and travel times and develop strategies to mitigate the impacts.

Incorporating source terms into the traffic flow equation can be challenging, as these terms are often complex and may depend on a wide range of factors, including driver behavior, road conditions, and weather conditions. However, advances in modeling techniques and computational methods have allowed incorporation of these factors into traffic flow models. There is growing interest in developing more realistic and accurate models of traffic dynamics that incorporate source terms.

## 2 Higher-order traffic flow model with source terms: Extended ARZ model

### 2.1 AW and Rascle model

AW and Rascle [2] proposed a new anisotropic continuum traffic flow model, addressing all the drawbacks of the Payne model [24] pointed out by Daganzo [4]. AWR introduced a convective derivative of pressure instead of a spatial derivative. This model considers not only the spatial changes of density but also the spatiotemporal changes in density. The model includes (1) and the velocity dynamics equation given in (3).

$$\frac{\partial(v+p(k))}{\partial t} + v\frac{\partial(v+p(k))}{\partial x} = \frac{V(k)-u}{T} \qquad (3)$$

Here $p(k)$ is the traffic density, an increasing function of density, expressed as (4).

$$p(k) = C^2 k^\gamma \qquad (4)$$

Here $C$ is a constant equal to 1. $V$ represents the preferred velocity proposed by [2], and $v$ is the vehicle traveling speed. In the case of traffic flow, $V(k)$, the preferred velocity is assumed to be a steady state speed corresponding to $k$ density. By choosing the pressure function as (4), AWR addressed the two most crucial criticisms of the Payne-type model, i.e., anisotropy and negative speeds at the tail of congested regions [4,15]. This model nicely predicts the instabilities near vacuum, i.e., for very light traffic. Since the predictive model ability largely depends on proper pressure function selection, Moutari and Rascle [20] formulated different pressure terms for the AWR model.

$$p(k) = \begin{cases} \dfrac{v_{ref}}{\gamma}\left(\dfrac{k}{k_{jam}}\right)^\gamma, & if\ \gamma > 0 \\[3mm] -v_{ref}\ln\left(\dfrac{k}{k_{jam}}\right), & if\ \gamma = 0 \end{cases}$$

Here $k$ denotes the fraction of space occupied by vehicles, $\gamma$ is a constant that can be equal to or greater than zero, $v_{ref}$ is a given reference velocity, which is equivalent to the maximum observed speed ($v_{max}$), and $k_{jam}$ is the maximal density as a fraction equal 1.

### 2.2 Incorporate source terms in continuum equation

Consider a road section with ramps to facilitate vehicle entry and exit from the main road of lengths $a$ to $b$. The vehicles may exit the road along the interval $[a,b]$. $g_{out}(t)$ is the fraction of traffic density per unit time that exits along $[a,b]$ and can be expressed as follows (5) [3]:

$$a_{out}(t,x) = g_{out}(t)\chi_{[a,b]}(x) \tag{5}$$

The sink term is negatively proportionate with the existing and jam densities, as traffic experiences exits from the road when fully crowded. Therefore, the sink term can be defined as the rate of exit flow per unit length of the road, introduced as follows (6):

$$sink = -\frac{a_{out}(t,x)k}{k_{jam}} \tag{6}$$

Similarly, if there is an entry along $[a,b]$, the entry flow will linearly decrease with density. Here $\chi_{[a,b]}(x)$ represents the localization of the ramp between $x = a$ and $x = b$ [3]. The source term can be defined as the rate of entry flow per unit length of the road as follows (7). $a_{in}$ and $a_{out}$ are the maximum incoming and outgoing flows during full free-flow and congested conditions:

$$source = a_{in}(t,x)\left(1 - \frac{k}{k_{jam}}\right)$$

$$a_{in}(t,x) = g_{in}(t)\chi_{[a,b]}(x) \tag{7}$$

Then, the model for entries/exits that have a convective part system is as follows (8).

$$\partial_t k + \partial_x(kv) = a_{in}(t,x)\left(1 - \frac{k}{k_{jam}}\right) - \frac{a_{out}(t,x)k}{k_{jam}} = \left(a_{in} - \frac{(a_{in}+a_{out})k}{k_{jam}}\right) \tag{8}$$

Let's assume the source term as $A(k) = \left(a_{in} - \frac{(a_{in}+a_{out})k}{k_{jam}}\right)$. This study assumes the entry and exit flow from the stream will depend on the ratio of existing density to jam density $\frac{k}{k_{jam}}$. This ratio represents the overall occupancy of the road in any instance. A ratio of $\frac{k}{k_{jam}}$ close to one indicates that the road is fully occupied, so most vehicles will attempt to exit without having a chance to enter. Similarly, $\frac{k}{k_{jam}}$ a ratio close to zero encourages entries and discourages exits.

### 2.3 Continuum & velocity dynamics equation

A relaxation term is added to ensure that the velocity of the flow is relaxed towards the preferred speed. So, the velocity dynamics equation is given by (9) (Payne type) [6,7].

$$\partial_t\big(v + p(k)\big) + v\partial_x\big(v + p(k)\big) = \frac{V(k)-v}{\delta} \qquad (9)$$

Here $V(k)$ is the maximum allowable speed, and $p(k)$ is the pressure term, a function of density. AW Rascle model expressed the pressure term as $p(k) = C_o^2 k^\gamma - \phi$ ( $C_o$ *is the sonic speed,* $\gamma, \phi$ *are constant terms*). To express (9) in the form of a conservation equation, the following mathematical operation is done: $\big((8) \times (v + p(k)), (9) \times k\big)$:

$$\big(v + p(k)\big)\partial_t k + \big(v + p(k)\big)\partial_x(kv) = \big(v + p(k)\big)A(k)$$

$$v\partial_t k + p\partial_t k + v^2\partial_x k + kv\partial_x v + pv\partial_x k + pk\partial_x k = \big(v + p(k)\big)A(k) \quad (10)$$

and

$$k\partial_t\big(v + p(k)\big) + vk\partial_x\big(v + p(k)\big) = k\frac{V(k)-v}{\delta} \qquad (11)$$

$$k(\partial_t v + p_k(k)\partial_t k) + kv(\partial_x v + p_k(k)\partial_x k) = k\,\frac{V(k)-v}{\delta}$$

After adding (10+11), we get (12)

$$(v\partial_t k + k\partial_t v + p\partial_t k + kp_k\partial_t k) + \big(v^2 + pv + kvp_k(k)\big)\partial_x k +$$

$$(2kv + pk)\partial_x v = \big(v + p(k)\big)A(k) + k\,\frac{V(k)-v}{\delta} \qquad (12)$$

$$\partial_t(k\big(v + p(k)\big)) + \partial_x(vk\big(v + p(k)\big)) = \big(v + p(k)\big)A(k) + k\,\frac{V(k)-v}{\delta}$$

Expressing (8) and (12) in the form of a generalized system of equations as (13):

$$\partial_t U + \partial_x F(U) = R(U) \qquad (13)$$

Where $U$, $F(u)$, and $R(U)$ are expressed as below:

$$U = \begin{bmatrix} k \\ k\big(v + p(k)\big) \end{bmatrix}, F(U) = \begin{bmatrix} kv \\ kv\big(v + p(k)\big) \end{bmatrix}, \qquad (14)$$

$$R(U) = \begin{bmatrix} A(k) \\ \big(v + p(k)\big)A(k) + k\,\frac{V(k)-v}{\delta} \end{bmatrix}$$

## 3       Numerical scheme

This section will discuss the numerical methods for solving the proposed higher-order model. The central dispersion (CD) of the MacCormack scheme, along with smoothing, is applied for spatial discretization. In traffic flow modeling, the second-order model typically takes the form of a non-linear partial differential equation that describes the evolution of the traffic density and velocity fields over time. The MacCormack scheme is a finite-difference method that can be used to discretize and solve this equation numerically. The MacCormack numerical scheme is essential for solving second-order traffic flow models because it provides a high-accuracy, stable, and efficient method for simulating traffic dynamics on roads and highways. The stability of the MacCormack scheme means that it produces a solution that does not blow up or oscillate with time, even in the presence of complex physical phenomena like turbulence or shock waves. The ability of the MacCormack scheme to handle both advection and diffusion terms in the equation is also crucial for many applications. Advection represents the transport of a quantity by a flow, while diffusion represents the spreading of a quantity due to turbulent mixing. The MacCormack scheme can handle both processes simultaneously, making it a powerful tool for modeling complex flows. Finally, the efficiency of the MacCormack scheme is important for traffic flow modeling because it allows for fast and computationally efficient simulations of large-scale traffic networks. This is important for applications such as traffic management, where real-time predictions and control strategies are required.

### 3.1     MacCormack Scheme

MacCormack Scheme is a classical discretization scheme for the numerical solution of hyperbolic partial differential equations [17]. Applying the MacCormack scheme to (13), we have the following predictor and corrector steps as (15):

$$Predictor: \overline{U}_i^{n+1} = U_i^n - \frac{\Delta t}{\Delta x}[F(U_{i+1}^n) - F(U_i^n)] + \Delta R(U_i^n)$$

$$Corrector: U_i^{n+1} = 0.5[\overline{U}_i^{n+1} + U_i^n] - 0.5\frac{\Delta t}{\Delta x}[\overline{F}_i^{n+1} - \overline{F}_{i-1}^{n+1}] \qquad (15)$$

**Central dispersion Method:** The first method for smoothing, the CD method, is based on the Lax Fredrichs scheme [26] and can be implemented on the updated results from the corrected step of MacCormack scheme (16):

$$U_{i(CD)}^{n+1} = (1-S)U_i^{n+1} + S\frac{(U_{i+1}^{n+1}+U_{i-1}^{n+1})}{2} \qquad (16)$$

Where $S$ represents the weight parameter related to smoothing, varies $0 \leq S \leq 1$, and $S = 0$ refers to no smoothing, reducing the method to the standard MacCormack.

### 3.2 Numerical stability

*CFL condition:* Numerical analysis of partial differential equations is required to satisfy the convergence condition proposed by Courant–Friedrichs–Lewy. In the case of finite-difference approximation, the time and space step of the numerical scheme must satisfy the Courant–Friedrichs–Lewy conditions to avoid numerical instabilities. The Courant–Friedrichs–Lewy condition for a one-dimensional case has to be followed the following form (17):

$$CFL = v\frac{\Delta t}{\Delta x} \le CFL_{max} \qquad (17)$$

where $CFL$ is the Courant number, $v$ is the magnitude of velocity, $\Delta t, \Delta x$ time steps, space step. $CFL_{max}$ is the maximum value of the Courant number, typically $CFL_{max} = 1$ for explicit solving or must be smaller than 1. Otherwise, the numerical viscosity would be negative. In the case of second-order traffic flow models, the speed used in the CFL number calculation is equal to the maximum of two Eigenvalues, which is greater than the backward propagation speed.

## 4 Numerical test

In this section, numerical tests are performed to evaluate the proposed numerical methods for application on Aw-Rascle models [2]. The model parameters take values given in Table 1. All test cases are performed on a homogeneous road of length L to confine the effects of additional source terms. At the boundaries, the following conditions (18) are applied:

$$U_1^n = U_2^n, \ and \ U_l^n = U_{l-1}^n \qquad (18)$$

This condition (18) states that the density at the first cell (index 1) equals the value at the adjacent second cell (index 2). The second condition indicates that the density at the last cell is equal to the density at the cell right before it. These boundary conditions ensure that there are no sudden variations or discontinuities in density at the boundaries of the simulated traffic segment. The initial conditions for starting the numerical solutions are given below.

$$IC1: k(x,0) = \begin{cases} 0.46k_{jam}, \ if \ x < \frac{l}{2} \\ 0.10k_{jam}, \ if \ x \ge \frac{l}{2} \end{cases} \qquad IC2: k(x,0) = \begin{cases} 0.90k_{jam}, \ if \ x < \frac{l}{2} \\ 0.55k_{jam}, \ if \ x \ge \frac{l}{2} \end{cases}$$

IC1 shows a sudden shift from congestion to free-flowing traffic as traffic density transitions from high to low levels. IC2 illustrates the transition from high to low densities within congested traffic settings. Both scenarios rely on specific initial traffic conditions that lead to the formation and propagation of expansion and shock waves. As a result, it becomes imperative to assess both the proposed model's effectiveness and the numerical scheme's capabilities in predicting traffic states.

**Table 1.** Physical and Numerical Parameters used for Aw-Rascle Models

| Parameter | Value |
|---|---|
| Length of the road (L) | 1200 m |
| Maximum speed ($v_{max}$) | 30 m/s |
| Maximum density ($k_{jam}$) | 0.15 veh/m |
| Critical density ($k_{cr}$) | $0.2667 k_{jam}$ veh/m |
| $C_0^2$ | $80 \ m^2/s^2$ |
| $\phi$ | $31.9 \ veh.m/s^2$ |
| $\gamma$ | 0.5 |
| $\Delta t$ | 1 sec |
| $\Delta x$ | 30 m |
| S | 0.01 |
| $a_{in}$ | 0.003 veh/m/sec |
| $a_{out}$ | 0.001 veh/m/sec |
| $\delta$ | 2 sec |
| C | 0.9 |

The behavior of the standard McCormack scheme for simulating Aw-Rascle-type models was investigated. We apply the method for some test cases to illustrate how the different source terms affect the solution.

**Case 1: no source and relaxation terms**

For this case $R(U) = 0$ in (13):

$$R(U) = \begin{bmatrix} 0 \\ 0 \end{bmatrix} \tag{19}$$

**Case 2: no source term but relaxation terms**

The $R(U)$ matrix for this case can be expressed as in (13):

$$R(U) = \begin{bmatrix} 0 \\ k\frac{V(k)-v}{\delta} \end{bmatrix} \qquad (20)$$

**Case 3: source term but no relaxation terms**

The $R(U)$ matrix for this case can be expressed as in (13):

$$R(U) = \begin{bmatrix} A(k) \\ (v + p(k))A(k) \end{bmatrix} \qquad (21)$$

**Case 4: with source term and relaxation terms**

The $R(U)$ matrix for this case can be expressed as in (13):

$$R(U) = \begin{bmatrix} A(k) \\ (v + p(k))A(k) + k\frac{V(k)-v}{\delta} \end{bmatrix} \qquad (22)$$

Figure 1-2 shows the test results for the modified AWR model with CD for two initial conditions. Both results show a lenient oscillation and smooth state transitions. Case-1,2 presents the model without source terms, whereas case3,4 includes the source terms in the model. Figures 1-2(a,b) show the traffic state propagation for IC1 and IC2 without incorporating the source terms. The model generates expansion and shock waves based on the initial traffic conditions, which support the KW model theory. Since this study applied the McCormack scheme with smoothing instead of distinct wave boundaries, Figure 1-2(a,b) shows smooth transitions from one state to another. Also, Figure 1-2(a,b) shows a dark color at x=6000 over time, which depicts the oscillation during state transition. In Figure 1-2(c,d), localized variations impacting specific points are simulated using source terms in time-space. As the source and sink terms are adjusted with the prevailing traffic state, one can see the smooth transition of dense to light or light to dense traffic in the time space. Figure 1-2(c,d) demonstrated less oscillation by reducing dark density at the state transition location. Figure 3-4 exhibits the cross-section view of Figure 1-2 at 100sec and 200sec. In case-4, oscillation is higher for both cases, which discourages the inclusion of source and relaxation terms in the modified model. For both IC, it is evident from Figure 3-4 that case 3 produces better state transition and mild oscillations than the other three cases. Results also show that the initial instabilities and oscillations decrease over time, especially in the cross-sectional plots for cases 3 and 4 (Figures 3 and 4). In the numerical tests, the parameter values are chosen as in Table 1, similarly as given in [19].

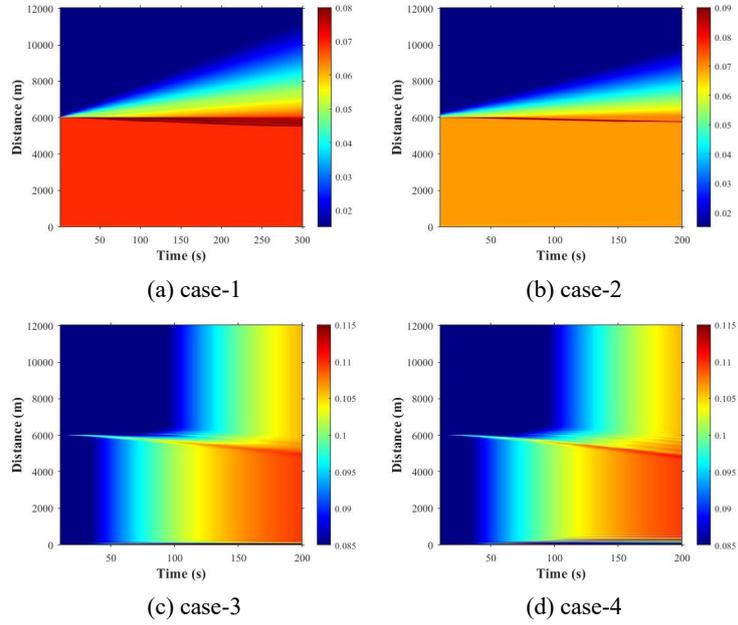

(a) case-1

(b) case-2

(c) case-3

(d) case-4

Fig.1: IC-1: Spatiotemporal evolution of local density transition from higher congested ($k > k_{cr}$) to lower free-flow density ($k < k_{cr}$) for AWR model after applying modified MacCormack-CD scheme

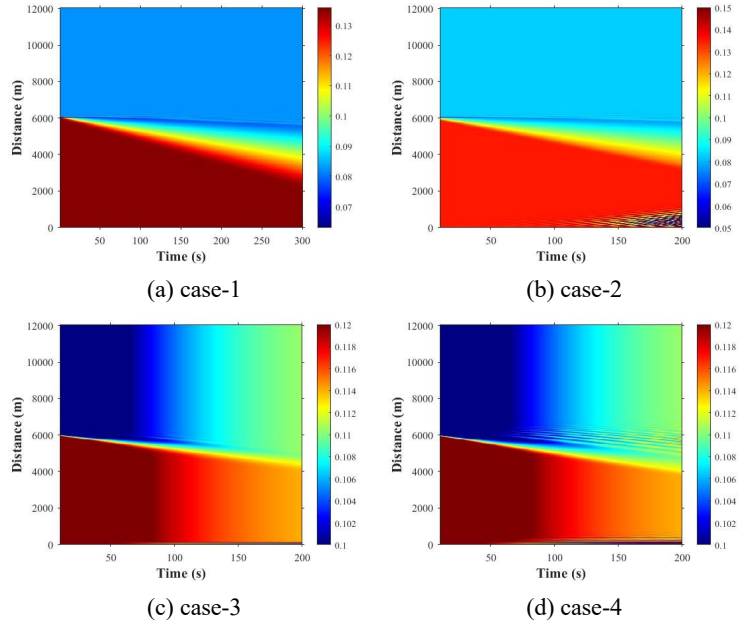

(a) case-1                 (b) case-2

(c) case-3                 (d) case-4

Fig.2: IC-2: Spatiotemporal evolution of congested traffic with a transition from high ($k > k_{cr}$) to low densities ($k > k_{cr}$) for AWR model after applying modified MacCormack-CD scheme.

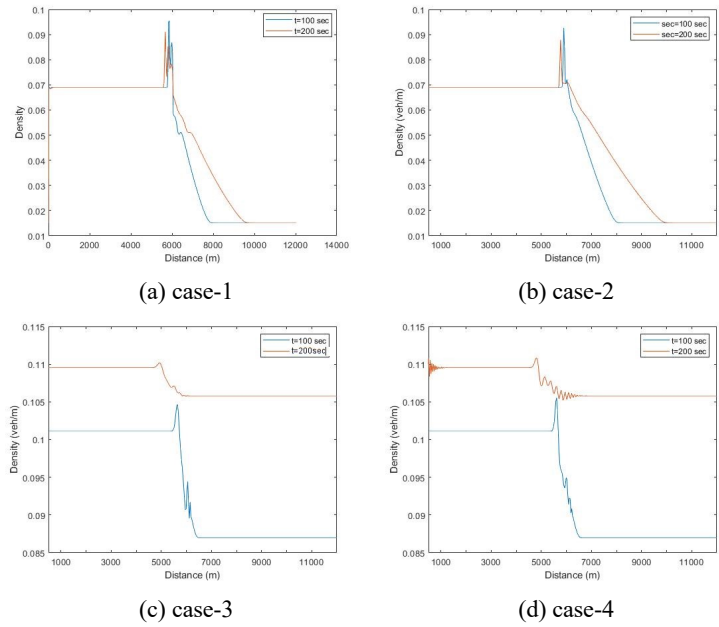

(a) case-1                    (b) case-2

(c) case-3                    (d) case-4

Fig.3: IC-1, congested traffic with a transition from high ($k > k_{cr}$) to low free flow density ($k < k_{cr}$) for AWR model for cross-section after applying modified MacCormack-CD scheme at $t = 100sec$ and $t = 200sec$.

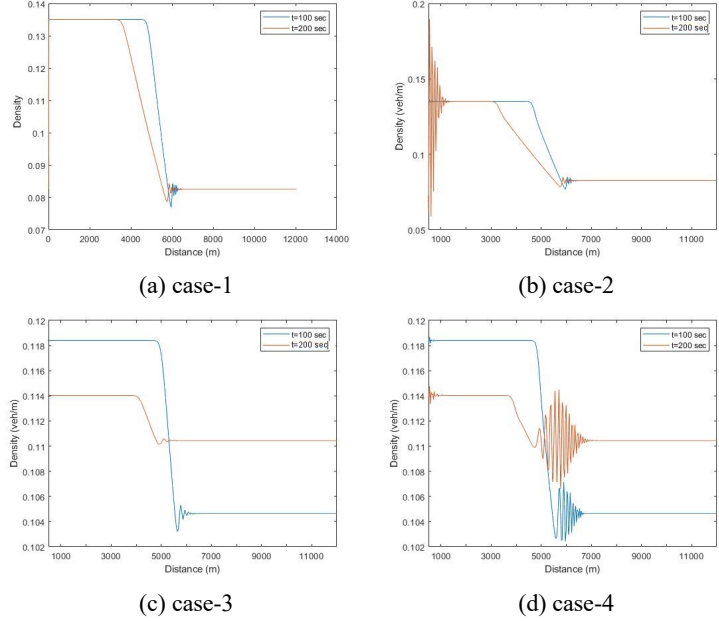

(a) case-1  (b) case-2

(c) case-3  (d) case-4

Fig.4: IC-2, congested traffic with a transition from high ($k > k_{cr}$) to low density ($k > k_{cr}$) for AWR model for cross-section after applying modified MacCormackCD scheme at $t = 100sec$ and $t = 200sec$.

## 5    Validation

The traffic states were estimated using Wi-Fi sensors installed as roadside units (RSU) and onboard units (OBU). More details about the data collection and traffic-state estimation can be found in the literature [18]. The Wi-Fi sensors were installed on the T-Nagar main road in Chennai. The travel time of the vehicles was calculated by matching the unique MAC id of Wi-Fi-enabled devices at different places passing through the vicinity of the sensors on the road. Wi-Fi sensors were installed at every 100m of the road. In this way, the speed of the stream can be determined from the Wi-Fi matching time and location. We calculated the density of the stream by considering the Greenshield model and assuming traffic moves in the equilibrium state. A comparison of empirical observations with model predictions is shown in Figure 5, where the colormap of density plots shows how well the prediction matches the empirical data. Figure 6 compares the model prediction accuracy with the empirical data sets at time $t = 10s$ and $t = 15s$. The discrepancy between predicted outcomes and actual observations arises because this study's simulated source and sink terms depend on cur-rent density conditions rather than the true field conditions. The prediction is more accurate when the traffic state propagates with fewer fluctuations.

The numerical scheme does not capture sudden changes in the traffic states since the method produces stable and smooth predictions. We compare the predicted density's root-mean-square errors (RMSE) with the empirically observed density in Table 2. The RMSE for data-1 and data-2 at 10 sec is 0.066 and 0.083 veh/m, respectively, increasing with time to 0.085 and 0.105 at 20 sec. Therefore, this model can predict the heterogeneous traffic state with minimal error, magnifying the error over time.

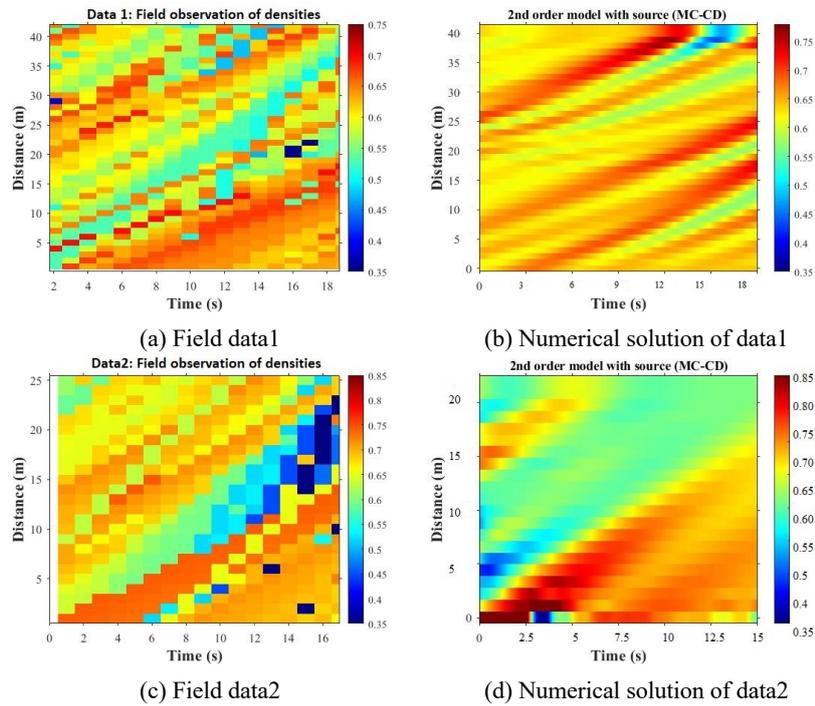

(a) Field data1

(b) Numerical solution of data1

(c) Field data2

(d) Numerical solution of data2

Fig.5: Spactiotemporal evolution of local density for the initial and boundary condition of data1 and data2 for extended AWR model after applying modified MacCormack-CD scheme. Here (a,b) represents field observations, and (c,d) shows prediction density.

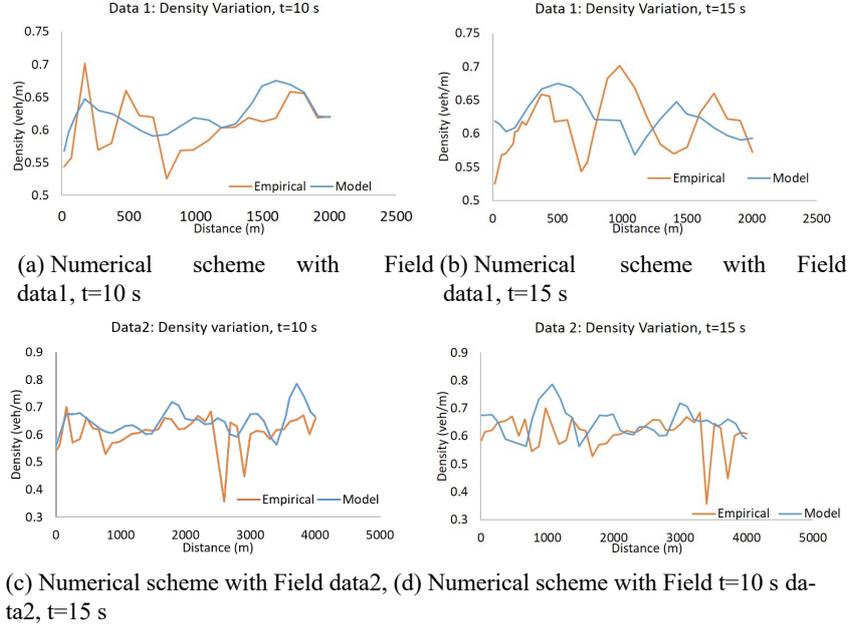

(a) Numerical scheme with Field data1, t=10 s

(b) Numerical scheme with Field data1, t=15 s

(c) Numerical scheme with Field data2, t=15 s

(d) Numerical scheme with Field t=10 s data2, t=15 s

Fig.6: Comparison of model performance prediction with the field density for data 1 and data 2 after applying MacCormack-CD scheme at t=10 s and t=15 s.

**Table 2.** Root Mean Square Error (RMSE) of empirical observation and the model prediction at 10 sec and 20 sec.

| Data | t = 10 s | t = 15 s |
|------|----------|----------|
| Data 1 | 0.066 | 0.085 |
| Data 2 | 0.083 | 0.105 |

## 6    Discussions

The study presents a new approach to modeling the traffic flow by incorporating source terms into the traditional continuum model. The source terms represent the external influences on the traffic flow, such as the inflow and outflow of vehicles from entry/exit ramps. This article proposes an in-homogeneous continuum equation with a velocity dynamics equation including diffusion term to model the link traffic flow, allowing entries and exits. While several studies have used higher-order traffic flow models to incorporate inertial effects, driver anticipation, and relaxation to a preferred speed only in the velocity dynamics, this article presented a second-order model with variable source and sink term in the continuum equation. The proposed model is numerically analyzed with the Mac-Cormack central dispersion smoothing scheme to prevent the numerical oscillation for spatial accuracy. The simulation test was performed for eight cases with or without source and relaxation terms for two initial conditions. Case-3, i.e., the model with only source terms, showed a smooth state transition without much oscillation and inflection. However, case-4, relaxation with source terms, predicts as similar to case-3 but was sometimes ignored due to its oscillating nature in state transition. In traffic flow modeling, relaxation terms are often used to account for the effects of driver behavior and vehicle interactions, allowing the simulated flow to approach a desired equilibrium speed. However, relaxation terms in the velocity dynamic equation can potentially lead to oscillations, depending on their specific formulation and the numerical methods used in the simulation. Improper or excessive use of relaxation terms can introduce numerical instability and oscillations in the simulation results. These oscillations can manifest as unrealistic fluctuations in vehicle velocities or densities. Hence future research is needed in the direction of identifying appropriate numerical methods and carefully calibrating and validating the model parameters. The method was validated with the empirical traffic density data collected from a road stretch. The comparison implies that the model can accurately predict the trend of traffic state propagation with a minimal error in the magnitude. This error in density magnitude prediction occurred due to the assumption of source and sink terms as a function of prevalent traffic density.

The main contributions of the paper are as follows:

- Incorporating source and sink terms: By adapting the AW Rascle-type velocity dynamic equation and modifying the continuum equation, the study significantly simplifies the system of equations. This integration enriches the model's ability to account for the impact of various localized events or disturbances that are frequently observed in mixed traffic conditions.
- Numerical testing with McCormack scheme: Incorporation of source terms plays a vital role in creating smoother state transitions within the traffic flow. This smoothness is achieved without inducing excessive oscillations or inflections, signifying the improved stability and reliability of the modified model.
- Validation with empirical data: To enhance the real-world relevance of the proposed model, the study validates its predictions using empirical mixed traffic data, which demonstrates the practical applicability and accuracy of the model in capturing the intricacies of actual traffic conditions.

Overall, the paper provides a valuable contribution to the traffic modeling and simulation field. The study's findings suggest incorporating source terms in the higher-order traffic flow models makes it more efficient to simulate various localized events or disturbances that are frequently observed in mixed traffic conditions.

# References

1. Abgrall, R., Qiu, J.: High order methods for cfd problems. Journal of computational-physics (Print) 230(11) (2011)
2. Aw, A., Rascle, M.: Resurrection of" second order" models of traffic flow. SIAMjournal on applied mathematics 60(3), 916–938 (2000)
3. Bagnerini, P., Colombo, R.M., Corli, A.: On the role of source terms in continuum traffic flow models. Mathematical and computer modelling 44(9-10), 917–930 (2006)
4. Daganzo, C.F.: Requiem for second-order fluid approximations of traffic flow.Transportation Research Part B: Methodological 29(4), 277–286 (1995)
5. Delis, A., Nikolos, I., Papageorgiou, M.: High-resolution numerical relaxation approximations to second-order macroscopic traffic flow models. Transportation Research Part C: Emerging Technologies 44, 318–349 (2014)
6. Goatin, P.: The aw–rascle vehicular traffic flow model with phase transitions. Mathematical and computer modelling 44(3-4), 287–303 (2006)
7. Goatin, P., Laurent-Brouty, N.: The zero relaxation limit for the aw–rascle–zhangtraffic flow model. Zeitschrift für angewandte Mathematik und Physik 70(1), 31 (2019)
8. Helbing, D., Treiber, M.: Numerical simulation of macroscopic traffic equations.Computing in Science & Engineering 1(5), 89–98 (1999)
9. Holland, E.N., Woods, A.W.: A continuum model for the dispersion of traffic ontwo-lane roads. Transportation Research Part B: Methodological 31(6), 473–485 (1997)
10. Hoogendoorn, S., Bovy, P.: Multiclass macroscopic traffic flow modelling: a multilane generalisation using gas-kinetic theory. In: 14th International Symposium on Transportation and Traffic TheoryTransportation Research Institute (1999)
11. Klar, A., Greenberg, J.M., Rascle, M.: Congestion on multilane highways. SIAMJournal on Applied Mathematics 63(3), 818–833 (2003)
12. Laval, J.A.: Some properties of a multi-lane extension of the kinematic wave model(2003)
13. LeVeque, R.J.: Finite volume methods for hyperbolic problems, vol. 31. Cambridgeuniversity press (2002)
14. Lighthill, M.J., Whitham, G.B.: On kinematic waves ii. a theory of traffic flowon long crowded roads. Proceedings of the Royal Society of London. Series A. Mathematical and Physical Sciences 229(1178), 317–345 (1955)
15. Liu, G., Lyrintzis, A.S., Michalopoulos, P.G.: Improved high-order model for freeway traffic flow. Transportation Research Record 1644(1), 37–46 (1998)
16. Logghe, S., Immers, L.H.: Multi-class kinematic wave theory of traffic flow. Transportation Research Part B: Methodological 42(6), 523–541 (2008)
17. MacCormack, R.W.: The effect of viscosity in hypervelocity impact cratering. Journal of spacecraft and rockets 40(5), 757–763 (2003)
18. Maiti, N., Chilukuri, B.R.: Estimation of local traffic conditions using wi-fi sensortechnology. Journal of Intelligent Transportation Systems pp. 1–18 (2023)
19. Mohammadian, S., van Wageningen-Kessels, F.: Improved numerical method forawrascle type continuum traffic flow models. Transportation Research Record 2672(20), 262–276 (2018)

20. Moutari, S., Rascle, M.: A hybrid lagrangian model based on the aw–rascle traf-ficflow model. SIAM Journal on Applied Mathematics 68(2), 413–436 (2007)
21. Munjal, P., Hsu, Y.S., Lawrence, R.: Analysis and validation of lane-drop effectson multi-lane freeways. Transportation Research/UK/ (1971)
22. Munjal, P., Pipes, L.A.: Propagation of on-ramp density perturbations on unidirec-tional two-and three-lane freeways. Transportation Research/UK/ (1971)
23. Ngoduy, D., Liu, R.: Multiclass first-order simulation model to explain non-lineartraffic phenomena. Physica A: Statistical Mechanics and its Applications 385(2), 667–682 (2007)
24. Payne, H.: Models of freeway traffic and control, simulation councils proc. Ser.: Mathematical Models of Public Systems 1(1)
25. Strikwerda, J.C.: Finite difference schemes and partial differential equations. SIAM(2004)
26. Thomas, J.W.: Numerical partial differential equations: finite difference methods,vol. 22. Springer Science & Business Media (2013)
27. Wong, G., Wong, S.: A multi-class traffic flow model – an extension of lwr model-with heterogeneous drivers. Transportation Research Part A: Policy and Practice 36(9), 827 – 841 (2002)